\documentclass[11pt]{article}
\usepackage[latin1]{inputenc}
\usepackage{amsmath,amsthm,amsfonts}
\usepackage[numbers]{natbib}
\usepackage[colorlinks,ps2pdf]{hyperref}
\usepackage[numbers]{natbib}
\usepackage{graphicx,psfrag}
\usepackage{enumerate}

\theoremstyle{plain}
\newtheorem{theorem}{Theorem}
\newtheorem{lemma}{Lemma}
\theoremstyle{definition}
\newtheorem{remark}{Remark}
\newtheorem{example}{Example}

\newcommand{\reff}[1]{(\ref{#1})}
\newcommand{\lest}{\le_{\rm{st}}}
\newcommand{\ord}{\preceq}
\newcommand{\ordst}{\ord_{\rm{st}}}
\newcommand{\Z}{\mathbb{Z}}
\newcommand{\R}{\mathbb{R}}
\newcommand{\E}{\operatorname{E}}

\newcommand{\tsum}{\textstyle\sum}
\newcommand{\lambdanet}{\lambda_{\rm net}}
\newcommand{\erl}{\operatorname{Erl}}
\newcommand{\aErl}{a_{\rm{Erl}}}
\newcommand{\bErl}{b_{\rm{Erl}}}
\newcommand{\mumin}{\mu^{\rm{min}}}
\newcommand{\mumax}{\mu^{\rm{max}}}
\newcommand{\muminval}{\mu_{\rm{min}}}
\newcommand{\mumaxval}{\mu_{\rm{max}}}
\newcommand{\Xmp}{X^{\rm{mp}}}
\newcommand{\amp}{a^{\rm{mp}}}
\newcommand{\Smp}{S^{\rm{mp}}}
\newcommand{\thetamp}{\theta^{\rm{mp}}}
\newcommand{\bmp}{b^{\rm{mp}}}
\newcommand{\Xmpagg}{\hat X^{\rm{mp}}}
\newcommand{\Smpagg}{\hat S^{\rm{mp}}}

\begin{document}

\title{Stochastic bounds for two-layer loss systems}
\author{
 M. Jonckheere\thanks{
 Eindhoven University of Technology,
 PO Box 513, 5600 MB Eindhoven, The Netherlands. \url{http://homepages.cwi.nl/\~jonckhee/}}
 \and
 L. Leskelä\thanks{
 Helsinki University of Technology,
 PO Box 1100, 02015 TKK, Finland. \url{http://math.tkk.fi/\~lleskela/}}
}
\date{February 15, 2008}
\maketitle

\begin{abstract}
 This paper studies multiclass loss systems with two layers of servers,
 where each server at the first layer is dedicated to a certain customer class,
 while the servers at the second layer can handle all customer classes.
 The routing of customers follows an overflow scheme, where arriving customers are
 preferentially directed to the first layer.
 Stochastic comparison and coupling techniques are developed for studying
 how the system is affected by packing of customers,
 altered service rates, and altered server configurations.
 This analysis leads to computationally fast upper and lower bounds for the performance
 of the system.
\end{abstract}


\


\section{Introduction}
\label{sec:introduction}

This paper studies multiclass loss systems with two layers of
servers, where each server at the first layer is dedicated to a
certain customer class, and the servers at the second layer can
handle all customer classes. Arriving customers are routed to vacant
servers in one of the layers, with preference given to the first
layer; or rejected otherwise. This policy is commonly referred to as
overflow routing.

Layered networks with overflow routing are commonly used in telecommunications
services, because different layers of service may increase the system capacity.
In \emph{wireless communication networks} for
instance, the servers at the first layer correspond to radio
channels dedicated to a small geographical area (microcell), and the
second layer represents available radio channels in a larger area
covering several microcells; in \emph{telephone call centers}, the
first layer consists of call agents trained to handling certain
types of phone calls, and the second layer represents call agents
who are cross-trained to deal with all types of calls.

The analysis of multilayer loss systems is challenging even under the
simplest statistical assumptions, because the distributions of the
overflow processes from the first layer are complex, and the direct
numerical computation of the stationary distribution
is unfeasible even for relatively small systems
(Louth, Mitzenmacher, and Kelly~\cite{louth1994}).
Hence, approximative methods are needed for performance analysis (see
Kelly~\cite{kelly1991} for a broad overview). Classical
approximation techniques such as the equivalent random method and
the Hayward--Fredericks method~\cite{wolff1989}, and the recently
introduced hyperexponential decomposition (Franx, Koole, and
Pot~\cite{franx2006}), are based on parametrically modeling the
overflow processes from the first layer by simpler processes.
These methods have been observed to produce good approximations for many choices
of system parameters. However, they may require considerable amounts
of computation, and it is not clear whether they remain accurate over
the full parameter range.

The goal of this paper is to approximate the system via upper and lower bounds
that are easy to compute numerically, and conservative in the sense that the
true performance remains between the bounds for all choices of system parameters.
To construct the upper bound, we modify the system by redirecting
customers from the second layer into the first layer as soon as
servers become vacant. This so-called maximum packing policy causes
the number of customers per class to have a product-form stationary distribution
(Everitt and Macfadyen~\cite{everitt1983}).
The lower bound is constructed by moving all servers from the second layer into the first,
this way reducing the system into a product of independent Erlang loss models.

The main tools for proving the validity of the bounds are (i) Massey's theorem~\cite{massey1987}
characterizing the comparability of two Markov jump processes; and (ii) stochastic coupling, where
versions of the processes describing the number of customers in the original and the reference
system are constructed in such a way that the difference of the two processes remains positive
with probability one. Coupling techniques have been successfully used by several authors for
deriving stochastic bounds for loss systems: Whitt~\cite{whitt1981} analyzed several single-class
queueing systems; Smith and Whitt~\cite{smith1981} studied the merging of two loss systems
together; Nain~\cite{nain1990} focused on multiclass single-layer loss systems; and Hordijk and
Ridder~\cite{hordijk1987} studied a special case of the two-layer loss system where the first
layer is fully dedicated to a single customer class. This paper extends some of the above results
to general multiclass two-layer loss systems, the main contribution being in showing that maximum
packing leads to upper bounds for the time-dependent and stationary distributions of the number of
customers in the system. In the special case where the first layer is fully dedicated to a single
customer class, this result improves the upper bound obtained by Hordijk and
Ridder~\cite{hordijk1987}.

The paper is organized as follows. Section~\ref{sec:modelDescription} introduces the model details
and notation. In Section~\ref{sec:preliminaryResult} we prove a preliminary comparison result that
is key to analyzing the monotonicity of the system. Section~\ref{sec:stochasticComparisons}
analyzes how the time-dependent distribution of the system is affected by maximum packing,
different server configurations, and altered service rates, and in
Section~\ref{sec:boundsOfTheStationaryDistribution} we carry out a similar analysis for the system
in steady state. Section~\ref{sec:conclusion} concludes the paper.

\section{Model description}
\label{sec:modelDescription}

\subsection{Two-layer loss system with overflow routing}
\label{sec:modelOriginal}

\begin{figure}[h]
 \psfrag{L1}{ $\lambda_1$}
 \psfrag{L2}{ $\lambda_2$}
 \psfrag{L3}{ $\lambda_3$}
 \psfrag{M1}{$m_1$}
 \psfrag{M2}{$m_2$}
 \psfrag{M3}{$m_3$}
 \psfrag{N}{$n$}
 \psfrag{La1}{ Layer $1$}
 \psfrag{La2}{ Layer $2$}
 \centering
 \includegraphics[height=50mm]{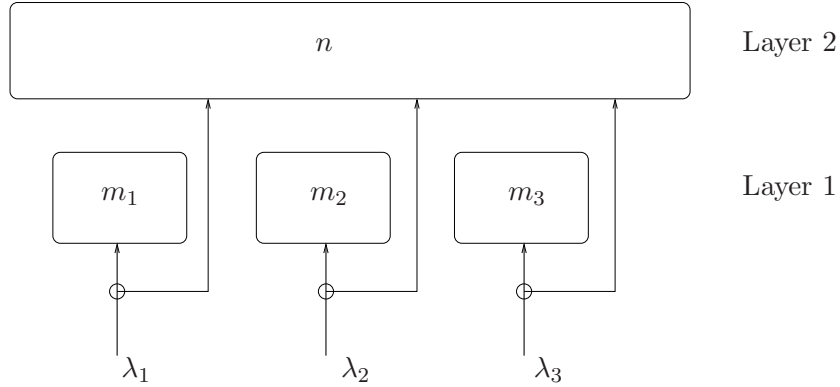}
 \caption{\label{fig:system} Two-layer loss network with three customer classes ($K=3$).}
\end{figure}

We consider a loss system with $K$ customer classes and two layers of servers, where layer~1
contains $m_k$ servers dedicated to class $k$, and layer~2 consists of $n$ servers capable of
serving all customer classes. Arriving class-$k$ customers are routed to vacant servers in one of
the layers, with preference given to layer~1; or rejected otherwise (Figure~\ref{fig:system}). For
analytical tractability, we assume that the interarrival times and the service requirements of
class-$k$ customers are exponentially distributed with parameters $\lambda_k$ and $\mu_k$,
respectively, and that all these random variables across all customer classes are independent. For
brevity, we denote $m=(m_1,\dots,m_K)$, $\lambda = (\lambda_1,\dots,\lambda_K)$, and $\mu =
(\mu_1,\dots,\mu_K)$.

Denote by $X_{i,k}(t)$ the number of class-$k$ customers being served at layer $i$ at time $t$.
The system is described by the continuous-time stochastic process $X = (X_{i,k})$ taking values in
\begin{equation}
 \label{eq:stateSpace}
 S=\{x \in \Z_+^K \times \Z_+^K: \ x_{1,k} \le m_k \ \forall k, \ \sum_{k=1}^K x_{2,k} \le n\}.
\end{equation}
Following the usual convention, we assume without loss of generality that all processes have paths
in the space $D(\R_+,S)$ of right-continuous functions with left-hand
limits~\cite{kallenberg2002}.

Let us denote by $e_{i,k}$ the unit vector in $\Z_+^K \times \Z_+^K$ corresponding to the
coordinate direction $(i,k)$. Moreover, define the sets
\begin{align}
 \label{eq:defA1}
 A_{1,k} &= \{x \in S: x_{1,k} < m_k\}, \\
 \label{eq:defA2}
 A_{2,k} &= \{x \in S: x_{1,k} = m_k, \ \sum_{l=1}^K x_{2,l} < n\}, \\
 \label{eq:defB}
 B_k     &= \{x \in S: x_{1,k} = m_k, \ \sum_{l=1}^K x_{2,l} = n\}.
\end{align}
The set $A_{i,k}$ represents the set of states where an arriving class-$k$ customer is assigned to
a layer-$i$ server, and $B_k$ is the set of states where arriving class-$k$ customers are
rejected. The process $X$ is a continuous-time Markov process on $S$ with the upward transitions
$x \mapsto x + e_{i,k}$ at rate $\lambda_{i,k}(x)$, and downward transitions $x \mapsto x -
e_{i,k}$ at rate $\phi_{i,k}(x)$, where
\begin{equation}
  \label{eq:generator}
  \begin{aligned}
    \lambda_{i,k}(x) &= \lambda_k 1(x \in A_{i,k}), \\
    \phi_{i,k}(x)    &= \mu_k x_{i,k}.
  \end{aligned}
\end{equation}

\subsection{Maximum packing}
\label{sec:modelMaximumPacking}

To approximate the original two-layer loss system, we consider a
modification of the system, where customers are redirected from
layer~2 to layer~1 as soon as servers become vacant. This
corresponds to the so-called maximum packing policy introduced by
Everitt and Macfadyen~\cite{everitt1983}. The process $\Xmp$
describing the number of customers in this system is a
continuous-time Markov process on $S$ with the upward transitions $x
\mapsto x + e_{i,k}$ at rate $\lambda'_{i,k}(x)$, and downward
transitions $x \mapsto x - e_{i,k}$ at rate $\phi'_{i,k}(x)$, where
\begin{equation}
  \label{eq:generatorPacking}
  \begin{aligned}
  \lambda'_{i,k}(x) &= \lambda_k 1(x \in A_{i,k}), \quad i=1,2, \\
  \phi'_{1,k}(x) &= \mu_k x_{1,k} 1(x_{2,k}=0), \\
  \phi'_{2,k}(x) &= \mu_k x_{1,k} 1(x_{2,k}>0) + \mu_k x_{2,k}.
\end{aligned}
\end{equation}

\begin{remark}
 \label{rem:maximumPacking}
 A remarkable property of the maximum packing policy is that
 all states outside the set
 $\Smp = \cap_{k=1}^K \{x \in S: x_{1,k}=m_k \ \text{or} \ x_{2,k}=0 \}$
 are transient for $\Xmp$. Moreover, note that for $x\in \Smp$,
 $x_{1,k} = m_k$ if and only if $ x_{1,k}+x_{2,k} \ge m_k$,
 which implies that
 \begin{equation}
  \label{eq:repacking1}
  \begin{aligned}
   x_{1,k} &= (x_{1,k} + x_{2,k}) \wedge m_k, \\
   x_{2,k} &= (x_{1,k} + x_{2,k} - m_k)^+.
  \end{aligned}
 \end{equation}
 As a consequence, the aggregate process $(\Xmp_{1,k} + \Xmp_{2,k})_{k=1}^K$
 tracking the total number of customers in each class,
 if started in $\Smp$, is equal in distribution to the Markov process
 $\Xmpagg=(\Xmpagg_1,\dots,\Xmpagg_K)$ on
 $\Smpagg = \{\hat x \in \Z_+^K: \sum_k (\hat x_k - m_k)^+ \le n\}$
 generated by the transitions
 \[
  \hat x \mapsto \left\{
  \begin{aligned}
   \hat x + e_k, &\quad \text{at rate} \ \lambda_k 1(\hat x + e_k \in \Smpagg), \\
   \hat x - e_k, &\quad \text{at rate} \ \phi_k(\hat x)=\mu_k \hat x_k.
  \end{aligned}
  \right.
 \]
 The structure of the above transition rates implies that the stationary distribution of $\Xmpagg$
 is a product of Poisson distributions truncated to $\Smpagg$~\cite{kelly1991}, which is easy to compute
 numerically. The stationary distribution of $\Xmp$ can then be recovered from that of $\Xmpagg$
 using the equalities~\reff{eq:repacking1}.
\end{remark}

\section{Preliminary result}
\label{sec:preliminaryResult}

This section establishes a general result that allows to compare two processes taking values in $S
\subset \Z_+^K \times \Z_+^K$ with respect to a specific preorder. This preorder, tailored to fit
the transition rates of the type in~\reff{eq:generator}, is defined by $x \ord y$, if $x_{1,k} \le
y_{1,k}$ for all $k$ and $|x| \le |y|$, where $|x| = \sum_{i,k} x_{i,k}$. For random variables
with values in $S$ we denote $U \ordst V$, if $\E \phi(U) \le \E \phi(V)$ for all bounded
measurable functions $\phi: S \to \R$ that are increasing with respect to the preorder $\ord$ on
$S$. Let us further extend these definitions to the Skorohod space $D(\R_+,S)$ of right-continuous
functions with left-hand limits by denoting $f \ord g$ if $f(t) \ord g(t)$ for all $t$. For
stochastic processes with paths in $D(\R_+,S)$ we denote $X \ordst Y$, if $\E \phi(X) \le \E
\phi(Y)$ for all bounded measurable maps $\phi: D(\R_+,S) \to \R$ that are increasing with respect
to the preorder $\ord$ on $D(\R_+,S)$. It will be clear from the context whether $\ord$ refers to
elements in $S$ or to functions in $D(\R_+,S)$.

Consider a continuous-time Markov process $X$ on $S \subset \Z_+^K \times \Z_+^K$ generated by the
transitions
\[
 x \mapsto \left\{
 \begin{aligned}
  x + e_{i,k} &\quad \text{at rate} \ \lambda_{i,k}(x), \\
  x - e_{i,k} &\quad \text{at rate} \ \phi_{i,k}(x),
 \end{aligned}
 \right.
\]
$i \in \{1,2\}$, $k \in \{1,\dots,K$\}, where $\lambda_{i,k}$ and $\phi_{i,k}$ are bounded nonnegative functions on $S$.
For consistency, we assume here that $\lambda_{i,k}(x)=0$ for all $x\in S$ such that $x+e_{i,k}\notin S$ and
$\phi_{i,k}(x)=0$ for all $x\in S$ such that $x-e_{i,k}\notin S$.
We assume that $Y$ is a similar process with state-dependent transition rates $\lambda'_{i,k}$ and $\phi'_{i,k}$.

\begin{theorem}
 \label{the:comparison}
 Let $X$ and $Y$ be Markov processes with paths in $D(\R_+,S)$ having
 upward transition rates $\lambda_{i,k}$ and $\lambda'_{i,k}$, and downward transition rates $\phi_{i,k}$ and $\phi'_{i,k}$,
 respectively. Assume that the following two conditions hold:
 \begin{enumerate}[(i)]
  \item For all $x,y\in S$ such that $x \ord y$ and $x_{1,k}=y_{1,k}$,
    \begin{align}
     \label{eq:comparisonMonoskillLambda}
     \lambda_{1,k}(x) &\le \lambda'_{1,k}(y), \\
     \label{eq:comparisonMonoskillPhi}
     \phi_{1,k}(x)    &\ge \phi'_{1,k}(y).
    \end{align}
  \item For all $x,y\in S$ such that $x \ord y$ and $|x| = |y|$,
    \begin{align}
     \label{eq:comparisonTotalLambda}
     \sum_{i,k} \lambda_{i,k}(x) &\le \sum_{i,k} \lambda'_{i,k}(y), \\
     \label{eq:comparisonTotalPhi}
     \sum_{i,k} \phi_{i,k}(x)    &\ge \sum_{i,k} \phi'_{i,k}(y).
    \end{align}
 \end{enumerate}
 Then $X \ordst Y$, given that the initial states satisfy $X(0) \ord Y(0)$.
\end{theorem}

\begin{proof}
 Denote the infinitesimal generators of $X$ and $Y$ by $p$ and $q$, respectively.
 Recall that $U \subset S$ is called an \emph{upper set}, if $x \in U$ and $x \ord y$ implies $y \in U$,
 and $V \subset S$ is called a \emph{lower set}, if the complement $V^c$ of $V$ is an upper set. Using
 a result of Massey~\cite[Theorem 5.3]{massey1987}\footnote{Massey
 formulated his result for partially ordered spaces, but all the proofs in his paper~\cite{massey1987}
 remain valid also for preorders that are not antisymmetric~\cite{last1995}.}
 (see also~\cite[Theorem 5]{kamae1977}), it suffices to verify that
 $p(x,U) \le q(y,U)$ for all $x \ord y$ and for all upper sets $U$ such that either $x \in U$ or $y \notin U$.
 Because $p(x,U) = -p(x,U^c)$ for all $x \in U$, this condition is equivalent to showing that for all $x \ord y$,
 \begin{equation}
  \label{eq:interToShow}
  \sum_{i,k} \lambda_{i,k}(x) \, 1(x+e_{i,k}\in U) \le \sum_{i,k} \lambda'_{i,k}(y) \, 1(y+e_{i,k}\in U)
 \end{equation}
 for all upper sets $U$ such that $x \notin U, y \notin U$, and
 \begin{equation}
  \label{eq:interToShowDecr}
  \sum_{i,k} \phi_{i,k}(x) \, 1(x+e_{i,k}\in V) \ge \sum_{i,k} \phi'_{i,k}(y) \, 1(y+e_{i,k}\in V)
 \end{equation}
 for all lower sets $V$ such that $x \notin V, y \notin V$.

 Assume $x \ord y$ and choose an upper set $U$ such that $x \notin U, y \notin U$.
 To verify the validity of~\reff{eq:interToShow}, let us consider separately the cases $|x| < |y|$ and $|x| = |y|$.
 Assume first $|x|<|y|$. Then $x + e_{1,k} \ord y$ for all $k$ such that $x_{1,k} < y_{1,k}$,
 and $x+e_{2,k} \ord y$ for all $k$.
 Hence because $U$ is an upper set and $y \notin U$, it follows that
 $x + e_{1,k} \in U$ only if $x_{1,k} = y_{1,k}$, and $x + e_{2,k} \notin U$ for all $k$. Thus,
 \begin{equation}
  \label{eq:restrictedSum}
  \sum_{i,k} \lambda_{i,k}(x) \, 1(x+e_{i,k}\in U)
  = \sum_{k: x_{1,k}=y_{1,k}} \lambda_{1,k}(x) \, 1(x+e_{1,k}\in U).
 \end{equation}
 Moreover, using inequality~\reff{eq:comparisonMonoskillLambda}, and noting that
 $x + e_{1,k} \ord y + e_{1,k}$ for all $k$ such that $y + e_{1,k} \in S$, we see that for
 all $k$ such that $x_{1,k} = y_{1,k}$,
 \begin{equation}
  \label{eq:lambdaLayer1}
  \lambda_{1,k}(x)  \, 1(x+e_{1,k} \in U) \le \lambda'_{1,k}(y) \, 1(y + e_{1,k} \in U).
 \end{equation}
 Substituting~\reff{eq:lambdaLayer1} into~\reff{eq:restrictedSum} shows the validity of~\reff{eq:interToShow}.

 Let us next focus on the case $|x| = |y|$. Note first that if
 $x + e_{1,l} \in U$ for some $l$ such that $x_{1,l} < y_{1,l}$, or
 $x + e_{2,l} \in U$ for some $l$, then
 $y + e_{i,k} \in U$ for all $i$ and $k$. Hence it follows that the right-hand side of~\reff{eq:interToShow}
 equals $\sum_{i,k} \lambda'_{i,k}(y)$,
 which in light of assumption~\reff{eq:comparisonTotalLambda} guarantees the validity of~\reff{eq:interToShow}.
 On the other hand, if $x + e_{2,k} \notin U$ for all $k$, and $x_{1,k} = y_{1,k}$ for all $k$ such that $x + e_{1,k} \in U$,
 then equation~\reff{eq:restrictedSum} holds.
 Assumption~\reff{eq:comparisonMonoskillLambda} again implies~\reff{eq:lambdaLayer1}, which together with~\reff{eq:restrictedSum} shows the validity of~\reff{eq:interToShow}.

 The proof is completed by carrying out an analogous reasoning for lower sets, which shows that
 assumptions~\reff{eq:comparisonMonoskillPhi} and~\reff{eq:comparisonTotalPhi} imply~\reff{eq:interToShowDecr}.
\end{proof}

\section{Pathwise stochastic comparison}
\label{sec:stochasticComparisons}

This section contains the main results for analyzing the time-dependent
distribution of the system. Assuming first that all service rates across different
customer classes are equal, we study how the system is affected by maximum packing (Section~\ref{sec:maximumPacking}) and different server configurations (Section~\ref{sec:differentServerConfigurations}).
Section~\ref{sec:monotonicityWithRespectToServiceRates} provides a monotonicity result
that allows to extend the analysis to the case where the service rates are not assumed
equal, and Section~\ref{sec:perClassBounds} describes bounds for the per-class number
of customers in the system.

Recall that the usual stochastic order~\cite{muller2002} between real random variables is defined
by denoting $U \lest V$, if $\E f(U) \le \E f(V)$ for all bounded measurable increasing real
functions $f$. Moreover, for stochastic processes with paths in the Skorohod space $D(\R_+,\R)$,
we denote $X \lest Y$ if $\E f(X) \le \E f(Y)$ for all bounded measurable functions $f: D(\R_+,\R)
\to \R$ that are increasing with respect to the natural pointwise order on $D(\R_+,\R)$. A
\emph{coupling} of two stochastic processes $X$ and $Y$ with paths in $D(\R_+,\R)$ is a stochastic
process $(\hat X, \hat Y)$ with paths in $D(\R_+,\R^2)$, having $X$ and $Y$ as its marginals.
Recall that by Strassen's theorem, $X \lest Y$ if and only if there exists a coupling $(\hat X,
\hat Y)$ of $X$ and $Y$ such that $\hat X(t) \le \hat Y(t)$ for all $t$ almost
surely~\cite{kamae1977}. Strassen's theorem can further be extended to processes with paths in
$D(\R_+,S)$, compared with respect to a given preorder~\cite{lindvall1999}.

\subsection{Maximum packing}
\label{sec:maximumPacking}

Let $X$ be the process describing the number of customers in the two-layer loss system defined in
Section~\ref{sec:modelOriginal}, and denote by $\Xmp$ the corresponding process for the maximum
packing policy defined in Section~\ref{sec:modelMaximumPacking}. Recall from
Section~\ref{sec:preliminaryResult} that the preorder $x \ord y$ is defined by $x_{1,k} \le
y_{1,k}$ for all $k$ and $|x| \le |y|$. The following theorem is the main result of the paper. It
allows to conclude that the stochastic processes $t \mapsto |X(t)|$ and $t \mapsto |\Xmp(t)|$
satisfy $|X| \lest |\Xmp|$, given that $X(0) \ord \Xmp(0)$.

\begin{theorem}
 \label{the:packing}
 Assume that all service rates $\mu_k$ are equal and that the initial states satisfy $X(0) \ord \Xmp(0)$.
 Then $X \ordst \Xmp$.
\end{theorem}

Example~\ref{exa:noSamplePathRepacking} below shows that a purely deterministic sample path
comparison is not sufficient for proving Theorem~\ref{the:packing}; hence probabilistic techniques
are needed. Example~\ref{exa:packing} in Section~\ref{sec:boundsOverall} further shows that the
statement of Theorem~\ref{the:packing} may not be true, if the service rates are not assumed
equal.

\begin{example}
 \label{exa:noSamplePathRepacking}
 Consider a two-class system ($K=2$) with one server at layer~1 assigned to class~1 ($m_1=1, m_2=0$)
 and one server at layer~2 ($n=1$). Denote by $X = (X_{i,k})$ a path of the process tracking the number of
 customers in the original two-layer loss system, and let $\Xmp$ be a corresponding sample path for the maximum packing
 policy. Assume that during the time interval $[0,6]$ there are four
 arriving customers each having service time equal to three:
 three  class-$1$ arrivals at time epochs 0, 2, and 4; and one class-$2$ arrival at time epoch~3.
 Given that both systems start empty, then $X(6) = e_{1,1}$ but $\Xmp(6) = 0$.
\end{example}

\begin{lemma}
 \label{the:monotoneLambda}
 The transition rates $\lambda_{i,k}(x)$ defined in~\reff{eq:generator}
 satisfy:
 \begin{enumerate}[(i)]
  \item For all $x \ord y$ and for all $k$ such that $x_{1,k} = y_{1,k}$,
  \begin{equation}
   \label{eq:lambda1}
   \lambda_{1,k}(x) \le \lambda_{1,k}(y).
  \end{equation}
  \item For all $x \ord y$ and for all $k$ such that $|x| = |y|$,
  \begin{equation}
   \label{eq:lambda2}
   \sum_{i,k} \lambda_{i,k}(x) \le \sum_{i,k} \lambda_{i,k}(y).
  \end{equation}
 \end{enumerate}
\end{lemma}
\begin{proof}
 The inequality~\reff{eq:lambda1} is clear, because $\lambda_{1,k}(x)$ only depends on $x_{1,k}$.
 Assume next that $x \ord y$ and $|x| = |y|$.
 Assume that $y \in B_k$ for some $k$, where $B_k$ is defined in~\reff{eq:defB}.
 Then $\sum_l y_{2,l} = n$, which implies that
 $\sum_l x_{2,l} = n$ and $x_{1,l} = y_{1,l}$ for all $l$. Thus $x \in B_k$. We may thus conclude that for all $k$,
 $ 1(x \notin B_k) \le 1(y \notin B_k)$.
 Hence it follows that
 \begin{align*}
  \sum_{i,k} \lambda_{i,k}(x)
  = \sum_k \lambda_k 1(x \notin B_k)
  \le \sum_k \lambda_k 1(y \notin B_k) = \sum_{i,k} \lambda_{i,k}(y),
 \end{align*}
 which shows the validity of~\reff{eq:lambda2}.
\end{proof}

\begin{proof}[Proof of Theorem~\ref{the:packing}]
 Let $\lambda_{i,k}(x)$ and $\phi_{i,k}(x)$ be the transition
 rates of $X$ as defined in~\reff{eq:generator}, and let
 $\lambda'_{i,k}(x)$ and $\phi'_{i,k}(x)$ be the corresponding rates for
 $\Xmp$ as defined in~\reff{eq:generatorPacking}.
 Because $\lambda_{i,k}'(x) = \lambda_{i,k}(x)$ for all $x$, the validity
 of~\reff{eq:comparisonMonoskillLambda} and~\reff{eq:comparisonTotalLambda}
 in Theorem~\ref{the:comparison} follow by Lemma~\ref{the:monotoneLambda}.
 For the downward transitions, note that for all $x \ord y$ such that $x_{1,k} = y_{1,k}$ for some $k$,
 $\phi_{1,k}(x) = \mu_1 x_{1,k} = \mu_1 y_{1,k} \ge \mu_1 y_{1,k} 1(y_{2,k}=0) = \phi_{1,k}'(y)$.
 Moreover, for all $x \ord y$ such that $|x| = |y|$,
 \[
  \sum_k (\phi_{1,k}(x) + \phi_{2,k}(x))
  = \mu_1 |x|
  = \mu_1 |y|
  = \sum_k (\phi_{1,k}'(y) + \phi_{2,k}'(y)),
 \]
 so conditions~\reff{eq:comparisonMonoskillPhi} and~\reff{eq:comparisonTotalPhi}
 of Theorem~\ref{the:comparison} are valid. Hence Theorem~\ref{the:comparison} yields the claim.
\end{proof}

\subsection{Different server configurations}
\label{sec:differentServerConfigurations}

This section studies the effect of moving one server from layer~1 to layer~2. As in
Section~\ref{sec:modelOriginal}, we denote by $X$ the process describing the number of customers
in the system with server configuration $m=(m_1,\dots,m_K)$ in layer~1, and $n$ servers in
layer~2. Let $Y$ by the process corresponding to the modified system where one class-$k$ server
from layer~1 has been replaced by a server in layer~2. We assume $k=1$ without loss of generality.
Let $m' = (m_1-1,m_2,\dots,m_K)$ and $n' = n+1$, and define the sets $S'$, $A_{1,k}'$ and $B_k'$
as in~\reff{eq:stateSpace}--\reff{eq:defB} with $m$ and $n$ replaced by $m'$ and $n'$,
respectively. Then $Y$ is a Markov process on $S'$ having transition rates of the
form~\reff{eq:generator} with $A_{i,k}$ replaced by $A'_{i,k}$.

Let us denote by $x_2 = \sum_k x_{2,k}$ the number of customers being served at layer~2. Assuming
that all service rates $\mu_k$ are equal, it follows that the process $(X_{1,1},\dots,X_{1,K};
X_2)$ is Markov. With a slight abuse of notation, we will redefine the state space by $S=
\{(x_{1,1}, \dots, x_{1,K}; x_2) \in \Z_+^K \times \Z_+: x_{1,k} \le m_k \ \text{for all} \ k, x_2
\le n\}$, and denote by $e_2$ the unit vector in $\Z_+^K \times \Z_+$ corresponding to the last
coordinate. We will redefine the sets $A_{i,k}, B_k, A'_{i,k}, B'_k$, and $S'$ in a similar way,
identifying $\sum_{k=1}^K x_{2,k}$ with $x_2$.

\begin{theorem}
 \label{the:training}
 Assume that all service rates $\mu_k$ are equal, and that the
 initial states satisfy $Y(0) - X(0) \in \Delta$, where
 $
  \Delta = \{ 0, \, e_2, \, e_2 - e_{1,1}, \, 2 e_2 - e_{1,1} \}.
 $
 Then the stochastic processes $ t \mapsto |X(t)|$ and $t \mapsto |Y(t)|$ satisfy
 $|X| \lest |Y|$.
\end{theorem}

\begin{proof}
 Because $|x| \le |y|$ for all $x\in S$ and $y\in S'$ such that $y-x \in \Delta$, it is sufficient to construct a coupling~\cite{thorisson2000} of $X$ and $Y$ that takes
 values in $S_\Delta = \{ (x,y) \in S \times S': y-x \in \Delta \}$.
 Let $(\tilde X, \tilde Y)$ be a continuous-time Markov process on $S_\Delta$ generated by the joint arrivals
 \begin{alignat}{3}
  \label{eq:trainingA1}
  (x,y) &\mapsto (x + e_{1,k}, \, y + e_{1,k}) &\quad & \text{at rate} & \ \ & \lambda_k 1(x \in A_{1,k}, \, y \in A'_{1,k}), \\
  \label{eq:trainingA2}
  (x,y) &\mapsto (x + e_{1,k}, \, y + e_2)     &&       \text{at rate} &&      \lambda_k 1(x \in A_{1,k}, \, y \in A'_{2,k}), \\
  \label{eq:trainingA3}
  (x,y) &\mapsto (x + e_{1,k}, \, y)           &&       \text{at rate} &&      \lambda_k 1(x \in A_{1,k}, \, y \in B'_k), \\
  \label{eq:trainingA4}
  (x,y) &\mapsto (x + e_2,     \, y + e_2)     &&       \text{at rate} &&      \tsum_l \lambda_l 1(x \in A_{2,l}, \, y \in A'_{2,l}), \\
  \label{eq:trainingA5}
  (x,y) &\mapsto (x + e_2,     \, y)           &&       \text{at rate} &&      \tsum_l \lambda_l 1(x \in A_{2,l}, \, y \in B'_l), \\
  \label{eq:trainingA6}
  (x,y) &\mapsto (x,           \, y + e_2)     &&       \text{at rate} &&      \tsum_l \lambda_l 1(x \in B_l,     \, y \in A'_{2,l}),
 \end{alignat}
 and joint departures
 \begin{alignat}{3}
  \label{eq:trainingD1}
  (x,y) &\mapsto (x - e_{1,k}, \, y - e_{1,k})  &\quad & \text{at rate} & \ \ & \mu_1 y_{1,k}, \\
  \label{eq:trainingD2}
  (x,y) &\mapsto (x - e_{1,1}, \, y - e_2)      &&       \text{at rate} &&      \mu_1 (x_{1,1} - y_{1,1}), \\
  \label{eq:trainingD3}
  (x,y) &\mapsto (x - e_2,     \, y - e_2)      &&       \text{at rate} &&      \mu_1 x_2, \\
  \label{eq:trainingD4}
  (x,y) &\mapsto (x,           \, y - e_2)      &&       \text{at rate} &&      \mu_1 (y_{1,1} + y_2 - x_{1,1} - x_2).
 \end{alignat}

 Observe that all transition rates above are nonnegative, because $y_{1,1} \le x_{1,1}$ and $y_{1,1} + y_2 \ge x_{1,1} + x_2$,
 whenever $y-x \in \Delta$. To ensure that the transitions define a generator of a Markov process on $S_\Delta$, we need to verify
 that $y'-x' \in \Delta$ for all transitions $(x,y) \mapsto (x',y')$, where $y-x \in \Delta$.
 This is obvious for transitions \reff{eq:trainingA1}, \reff{eq:trainingA4}, \reff{eq:trainingD1}, and \reff{eq:trainingD3},
 because in these cases $y'-x' = y-x$. Let us consider the remaining cases one-by-one:
 \begin{itemize}
  \item If transition~\reff{eq:trainingA2} occurs, then $k=1$, because $y_{1,k} = x_{1,k}$ for all $k\neq 1$.
  Then $x_{1,1} < m_1$ and $y_{1,1} = m_1 - 1$, so it follows that either $y - x = 0$ or $y - x = e_2$. In both
  cases, $y' - x' \in \Delta$.
  \item If transition~\reff{eq:trainingA3} occurs, then again $k=1$. Then $x_{1,1} < m_1$ and $y_{1,1} = m_1 - 1$, which
  implies $y_{1,1} = x_{1,1}$. Moreover, $y_2 = n+1$, which is only possible if $y_2 = x_2 + 1$. Hence $y-x = e_2$, so
  that $y' - x' = e_2 - e_{1,1} \in \Delta$.
  \item If transition~\reff{eq:trainingA5} occurs, then $x \in A_{2,l}$ and $y \in B'_l$ for some $l$. Then
  $x_2 < n$ and $y_2 = n+1$, which implies that $y-x = 2 e_2 - e_{1,1}$. Hence $y' - x' = e_2 - e_{1,1} \in \Delta$.
  \item If transition~\reff{eq:trainingA6} occurs, then $x \in B_l$ and $y \in A'_{2,l}$ for some $l$.
  Then $x_2 = n$ and $y_2 < n+1$, so it follows that $y_2 = x_2$. Hence $y-x = 0$, and thus $y' - x' = e_2 \in \Delta$.
  \item If transition~\reff{eq:trainingD2} occurs, then $y_{1,1} < x_{1,1}$. Because $y - x \in \Delta$,
  this implies that either $y - x = e_2 - e_{1,1}$, so that $y' - x' = 0$; or $y - x = 2 e_2 - e_{1,1}$, so that
  $y' - x' = e_2$.
  \item If transition~\reff{eq:trainingD4} occurs, then $y_{1,1} + y_2 - x_{1,1} - x_2 > 0$. Because $y-x \in \Delta$, it follows that
  either $y - x = e_2$, so that $y' - x' = 0$; or $y - x = 2 e_2 - e_{1,1}$, so $y' - x' = e_2 - e_{1,1}$.
 \end{itemize}
 Hence, all transitions map $S_\Delta$ into $S_\Delta$, and the process $(\tilde X, \tilde Y)$ is well-defined.

 To show that $(\tilde X, \tilde Y)$ is a coupling of $X$ and $Y$, we must verify that
 the marginal transition rates of $(\tilde X,\tilde Y)$ match with the transition rates of $X$ and $Y$.
 Note first that the sum of transition rates such that $x \mapsto x + e_{1,k}$ is equal to $\lambda_k 1(x \in A_{1,k})$.
 Next, observe that $x \in A_{2,l}$ and $y-x \in \Delta$ imply that $y \notin A_{1,l}'$. Hence the sum of transition rates
 where $x \mapsto x + e_2$ is equal to
 \[
  \sum_l \lambda_l 1(x \in A_{2,l}, y \in A'_{2,l} \cup B'_l) = \sum_l \lambda_l 1(x \in A_{2,l}).
 \]
 Further, because the sum of all transition rates such that $x \mapsto x - e_{1,k}$ equals $\mu_1 x_{1,k}$ for all $k$, and the corresponding sum for $x \mapsto x - e_2$ is equal to $\mu_1 x_2$, we may conclude that the
 transitions of $\tilde X$ and $X$ occur at the same rates.

 Turning the attention to the rates of $\tilde Y$,
 note that $y-x \in \Delta$ and $y \in A_{1,1}'$ imply that $y_{1,1} < m_1 - 1$ and $x_{1,1} \le y_{1,1} + 1$,
 so it follows that $x \in A_{1,1}$. Moreover, $y-x \in \Delta$ and $y \in A_{1,k}'$ for $k\neq 1$ imply that
 $x_{1,k} = y_{1,k} < m_k$, so $x \in A_{1,k}$. Hence the total rate of transitions where
 $y \mapsto y + e_{1,k}$ is equal to $\lambda_k 1(x \in A_{1,k}, \, y \in A'_{1,k}) = \lambda_k 1(y \in A'_{1,k})$.
 Further, because the net rate of transitions where $y \mapsto y + e_2$ is equal to
 $\sum_l \lambda_l 1(y \in A'_{2,l})$, and because the corresponding net rates for $y \mapsto y - e_{1,k}$ and
 $y \mapsto y - e_2$ are equal to $\mu_1 y_{1,k}$ and $\mu_1 y_2$, respectively, we conclude that
 the transitions of $\tilde Y$ and $Y$ occur at the same rates. Hence, the process $(\tilde X, \tilde Y)$
 is a coupling of $X$ and $Y$.
\end{proof}

\subsection{Monotonicity with respect to service rates}
\label{sec:monotonicityWithRespectToServiceRates}

The results in Sections~\ref{sec:maximumPacking} and~\ref{sec:differentServerConfigurations} were
proved under the assumption that all service rates are equal. The following theorem describes a
monotonicity property that allows to compare systems not satisfying this assumption. Denote by $X$
the number of customers of the two-layer loss system defined in Section~\ref{sec:modelOriginal}.
Recall that the preorder $x \ord y$ is defined by $x_{1,k} \le y_{1,k}$ for all $k$ and $|x| \le
|y|$.

\begin{theorem}
 \label{the:monotonemu}
 Let $X^-$ and $X^+$ be modifications of the system with all service rates set to $\mumaxval = \max
 \mu_k$ and $\muminval = \min \mu_k$, respectively. Assume that the initial states satisfy
 $X^-(0) \ord X(0) \ord X^+(0)$. Then
 \[
  X^- \ordst X \ordst X^+.
 \]
\end{theorem}

\begin{remark}
 A simpler comparison statement, such as $|X| \lest |X^+|$
 given that $|X(0)| \le |X^+(0)|$, is not true in general. Using Massey's~\cite{massey1987} criteria for the
 preorder $|x| \le |y|$, it is not hard to check that
 a necessary condition for the above property
 is that $\sum_{i,k} \lambda_{i,k}(x) = \sum_{i,k} \lambda_{i,k}(y)$ whenever $|x| = |y|$. This equality
 fails for $x = \sum_k m_k e_{1,k} + (n-1)e_{2,1}$ and $y = x - e_{1,1} + e_{2,1}$.
\end{remark}

\begin{proof}[Proof of Theorem~\ref{the:monotonemu}]
 Note that $X^+$ has the same upward transitions as $X$ and downward transitions
 $\phi'_{1,k}(x) = \muminval x_{1,k}$, and  $\phi'_{2,k}(x) = \muminval x_{2,k}$.
 Now for all $x \ord y$ such that $x_{1,k} = y_{1,k}$ for some $k$,
 $\mu_k x_{1,k}  \ge \muminval x_{1,k}  = \muminval y_{1,k}$,
 and for all $x \ord y$ such that $|x| = |y|$,
 \begin{align*}
  \sum_k \mu_k ( x_{1,k} + x_{2,k} )
  &\ge \muminval \sum_k ( x_{1,k} + x_{2,k} ) = \muminval \sum_{i,k} ( y_{1,k} + y_{2,k} ),
 \end{align*}
 so conditions~\reff{eq:comparisonMonoskillPhi} and~\reff{eq:comparisonTotalPhi} of
 Theorem~\ref{the:comparison} are valid. Moreover, \reff{eq:comparisonMonoskillLambda} and~\reff{eq:comparisonTotalLambda} hold by Lemma~\ref{the:monotoneLambda},
 so Theorem~\ref{the:comparison} yields the claim for $X^+$. The claim for $X^-$
 is proved in a similar way.
\end{proof}

\subsection{Per-class bounds}
\label{sec:perClassBounds}

In this section, we prove upper and lower bounds for the per-class number of customers in the
system. Let $Z^s_{\lambda,\mu}(t)$ be the number of customers in the standard $s$-server Erlang
loss system at time $t$, defined as the right-continuous Markov process on $\{0,1,\dots,s\}$
having the upward transitions $x \mapsto x + 1$ at rate $\lambda 1(x<s)$ and the downward
transitions $x \mapsto x - 1$ at rate $\mu x$.

\begin{theorem}
 \label{the:perClassLower}
 Assume $Z^{m_k}_{\lambda_k, \mu_k}(0) = X_{1,k}(0)$. Then
 \begin{equation}
  \label{eq:perClassLower}
  Z^{m_k}_{\lambda_k, \mu_k} \lest X_{1,k} + X_{2,k}.
 \end{equation}
\end{theorem}
\begin{proof}
 Observe that the process $X_{1,k}$ tracking the number of class-$k$ customers being served at
 layer~1 has the same dynamics as a standard $m_k$-server Erlang loss system with arrival rate
 $\lambda_k$ and service rate $\mu_k$. Hence given $Z^{m_k}_{\lambda_k, \mu_k}(0) = X_{1,k}(0)$,
 the processes $Z^{m_k}_{\lambda_k, \mu_k}$ and $X_{1,k}$ have the same distribution, which
 immediately implies~\reff{eq:perClassLower}.
\end{proof}

\begin{theorem}
 \label{the:perClassUpper}
 Assume $X_{1,k}(0) + X_{2,k}(0) \le Z^{m_k+n}_{\lambda_k, \mu_k}(0)$. Then
 \begin{equation}
  \label{eq:perClassUpper}
  X_{1,k} + X_{2,k} \lest Z^{m_k+n}_{\lambda_k, \mu_k}.
 \end{equation}
\end{theorem}
\begin{proof}
 Assume without loss of generality that $k=1$. Let us construct a Markov process $(\tilde X, \tilde Y)$ on
 \[
  S_2 = \{ (x,y) \in S \times \{0,\dots,m_1+n\}: x_{1,1}+x_{2,1} \le y\}
 \]
 via the class-1 transitions for $i = 1,2$,
 \begin{alignat}{3}
  \label{eq:perClass1b}
  (x,y) &\mapsto (x + e_{i,1}, \, y+1) &\quad& \text{at rate} &\ \ & \lambda_1 1(x \in A_{i,1}, \, y < m_1 + n), \\
  \label{eq:perClass2b}
  (x,y) &\mapsto (x + e_{i,1}, \, y)   &&      \text{at rate} &&     \lambda_1 1(x \in A_{i,1}, \, y = m_1 + n), \\
  \label{eq:perClass3b}
  (x,y) &\mapsto (x,           \, y+1) &&      \text{at rate} &&     \lambda_1 1(x \in B_1,     \, y < m_1 + n), \\
  \label{eq:perClass4b}
  (x,y) &\mapsto (x - e_{i,1}, \, y-1) &&      \text{at rate} &&     \mu_1 x_{i,1}, \\
  \label{eq:perClass5b}
  (x,y) &\mapsto (x,           \, y-1) &&      \text{at rate} &&     \mu_1 (y - x_{1,1} - x_{2,1}),
 \end{alignat}
 and the class-$k$ transitions for $k\neq 1$ and $i = 1,2$,
 \begin{alignat}{3}
  \label{eq:perClass6b}
  (x,y) &\mapsto (x + e_{i,k}, \, y)   &\quad& \text{at rate} &\ \ & \lambda_k 1(x \in A_{i,k}), \\
  \label{eq:perClass7b}
  (x,y) &\mapsto (x - e_{i,k}, \, y)   &&      \text{at rate} &&     \mu_k x_{i,k}.
 \end{alignat}
 Note that all transition rates in~\reff{eq:perClass1b}~--~\reff{eq:perClass7b} are nonnegative for all $(x,y) \in S_2$.

 Let us now verify that all transitions map $S_2$ into $S_2$. Observe first that transition~\reff{eq:perClass2b} occurs only
 if $y=m_1+n$ and either $x_{1,1} < m_1$ or $\sum_{k=1}^K x_{2,k} < n$, which implies that $(x + e_{i,1}, \, y) \in S_2$ for $i=1,2$.
 Moreover, transition~\reff{eq:perClass5b} occurs only if $x_{1,1} + x_{2,1} < y$, so that $(x, y-1) \in S_2$. Clearly,
 all other transitions map $S_2$ into $S_2$. Thus the Markov process $(\tilde X, \tilde Y)$ on $S_2$ is well-defined.

 Moreover, the total rates of transitions in~\reff{eq:perClass1b}~--~\reff{eq:perClass7b} where
 $x\mapsto x + e_{i,k}$ and $x\mapsto x - e_{i,k}$ are equal to $\lambda_k 1(x \in A_{i,k})$
 and $\mu_k x_{i,k}$, respectively, for all $i$ and $k$. The corresponding total rates for
 $y \mapsto y+1$ and $y\mapsto y-1$ are equal to $\lambda_1 1(y < m_1 + n)$ and $\mu_1 y$, respectively.
 This shows that $(\tilde X, \tilde Y)$ is a coupling of $X$ and $Z_{\lambda_1,\mu_1}^{m_1 + n}$,
 so the inequality~\reff{eq:perClassUpper} holds.
\end{proof}

\begin{remark}
 \label{rem:perClassLower}
 The proof of Theorem~\ref{the:perClassLower} actually shows that inequality~\reff{eq:perClassLower}
 can be extended to arbitrary (random or nonrandom) arrival processes and service times.
 Example~\ref{exa:noSamplePathPerClass} shows why this kind of purely deterministic sample path comparison is not possible
 for obtaining the result in Theorem~\ref{the:perClassUpper}.
\end{remark}

\begin{example}
 \label{exa:noSamplePathPerClass}
 Consider a two-class system ($K=2$) with no servers at layer~1 ($m_1 = 0, m_2=0$)
 and one server at layer~2 ($n=1$). Denote by $X = (X_{i,k})$ a path of the process tracking the number of
 customers in the original two-layer loss system, and let $Z$ be a corresponding sample path of
 the modified (one-class) system that only accepts class-$1$ customers. Assume that during the time interval $[0,3]$ there are three
 arriving customers each having service time equal to two:
 a class-$2$ arrival at time epoch 0, and two class-$1$ arrivals at time epochs 1 and 2.
 Given that both systems start empty, then $X_{2,1}(3) = 1$ but $Z(3) = 0$.
\end{example}

\section{Bounds of the steady-state performance}
\label{sec:boundsOfTheStationaryDistribution}
Assume from now on that all arrival rates and service
rates are strictly positive, which implies that all Markov processes treated in the sequel have a
unique stationary distribution. In this section $\bar X = (\bar X_{i,k})$ denotes a random vector
describing the stationary number of class-$k$ customers being served at layer $i$ in the system,
and the quadruple $(m,n,\lambda,\mu)$ indicates that a performance quantity corresponds to a
system with server configuration $m=(m_1,\dots,m_K)$ at layer~1, $n$ servers at layer~2, arrival
rates $\lambda = (\lambda_1,\dots,\lambda_K)$, and service rates $\mu = (\mu_1,\dots,\mu_K)$.

\subsection{Per-class performance}
\label{sec:boundsPerClass}

Denote by $a_k = \E( \bar X_{1,k} + \bar X_{2,k})$ the stationary mean number of class-$k$
customers in the system, by $\theta_k$ the stationary mean class-$k$ throughput (the number of
class-$k$ customers completing service per unit time), and by $b_k$ the class-$k$ blocking
probability. Note that $a_k$ can be viewed as the mean class-$k$ work throughput (amount of
class-$k$ work served per unit time).

Let $\erl(s,\rho)$ be a random variable on $\{0,1,\dots,s\}$ having distribution $(\sum_{j=0}^s
\frac{\rho^j}{j!})^{-1} \frac{\rho^i}{i!}$, and denote its mean by $\aErl(s,\rho)$, and the
probability of being equal to $s$ by $\bErl(s,\rho)$. Note that $\bErl(s,\rho)$ is equal to the
famous Erlang~B formula.

\begin{theorem}
 \label{the:perClassSteadyState}
 The stationary number of class-$k$ customers in the system satisfies
 \begin{equation}
  \label{eq:perClassQueueLength}
  \erl(m_k, \lambda_k/\mu_k) \lest \bar X_{1,k} + \bar X_{2,k} \lest \erl(m_k + n, \lambda_k/\mu_k).
 \end{equation}
 Especially, the stationary class-$k$ mean number of customers
 is bounded by
 \begin{equation}
  \label{eq:perClassWorkThroughput}
  \aErl(m_k, \lambda_k/\mu_k) \le a_k \le \aErl(m_k + n, \lambda_k/\mu_k),
 \end{equation}
 the mean throughput by
 \begin{equation}
  \label{eq:perClassThroughput}
  \mu_k \aErl(m_k, \lambda_k/\mu_k) \le \theta_k \le \mu_k \aErl(m_k + n, \lambda_k/\mu_k),
 \end{equation}
 and the blocking probability by
 \begin{equation}
  \label{eq:perClassBlocking}
  \bErl(m_k+n,\lambda_k/\mu_k) \le b_k \le \bErl(m_k, \lambda_k/\mu_k).
 \end{equation}
\end{theorem}
\begin{proof}
 Let us consider a version of the process $X$ started at $X(0)= 0$,
 and let $Z^{m_k}_{\lambda_k, \mu_k}$ be as in Theorem~\ref{the:perClassLower} and $Z^{m_k+n}_{\lambda_k, \mu_k}$ be as in
 Theorem~\ref{the:perClassUpper}, both started at zero. Because all these
 processes are irreducible and positive recurrent, and because stochastic ordering is closed
 with respect to convergence in distribution~\cite{kamae1977}, the inequalities~\reff{eq:perClassQueueLength} follow
 by taking $t \to \infty$ in~\reff{eq:perClassLower} and~\reff{eq:perClassUpper}.

 The inequalities~\reff{eq:perClassWorkThroughput} follow by taking expectations, and the bounds~\reff{eq:perClassThroughput} are
 a consequence of $\theta_k = \mu_k a_k$.
 In light of the conservation laws $\lambda_k (1-b_k) = \theta_k$ and $\lambda_k (1-\bErl) = \mu_k \aErl$,
 these bounds in turn imply~\reff{eq:perClassBlocking}.
\end{proof}

\begin{figure}[h]
 \centering
 \includegraphics[width=.49\textwidth]{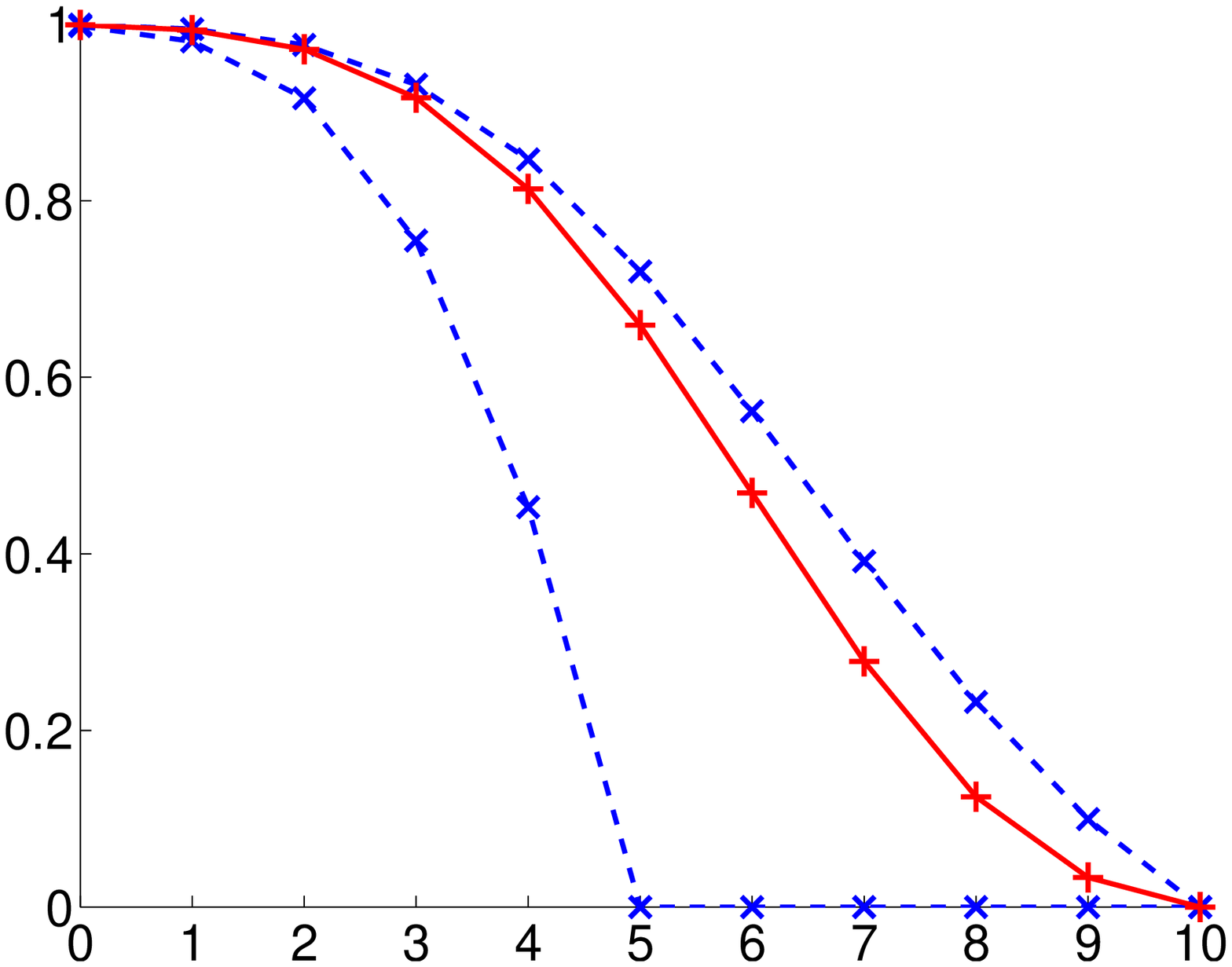}
 \includegraphics[width=.49\textwidth]{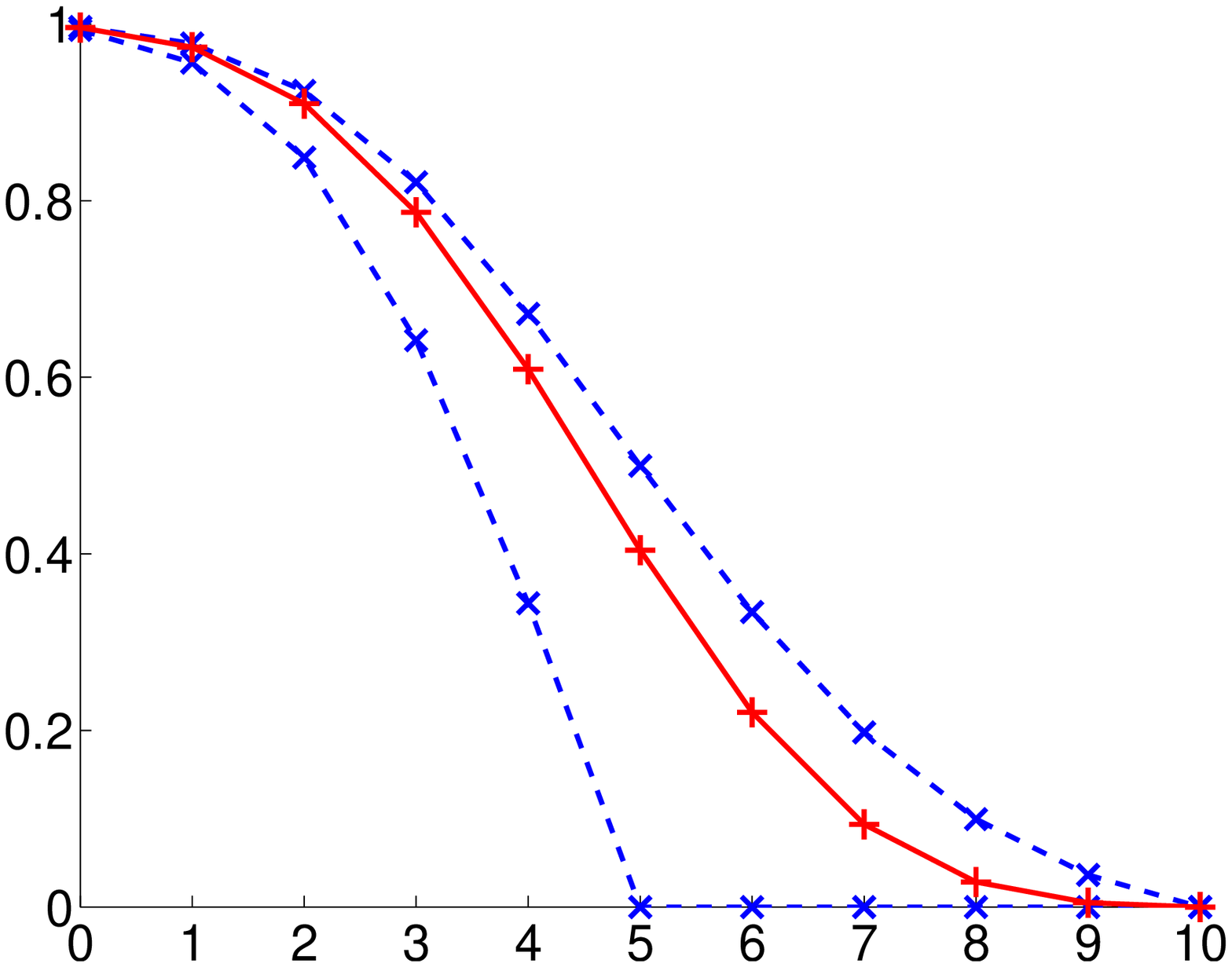}
 \caption{\label{fig:distributionNumberOfCustomers} Complementary cumulative distribution functions of the
 stationary number of class-1 customers (left) and class-2 customers (right), plotted for
 the original system (solid line) and for the Erlang bounds in~\reff{eq:perClassQueueLength} (dotted lines).}
\end{figure}

Figure~\ref{fig:distributionNumberOfCustomers} illustrates the bounds
in~\reff{eq:perClassQueueLength} for a loss network with server configuration $m=(5,5)$ and $n=5$,
where $\lambda = (7.5,7.5)$ and $\mu = (1, 1.3)$.

\begin{remark}
 The Erlang bounds~\reff{eq:perClassBlocking} for the blocking probability are well-known in the literature
 (see for example~\cite{borst2000}). The stochastic inequalities~\reff{eq:perClassQueueLength} can be viewed as
 extensions of these classical bounds.
\end{remark}

\subsection{Overall performance}
\label{sec:boundsOverall}

Denote by $a = \E |\bar X|$ the stationary mean total number of customers, by $\theta = \sum_k
\theta_k$ the stationary mean throughput, and by $b$ the stationary overall blocking probability.
Note that $a$ may be viewed as the mean work throughput (net amount of work served by the system
in unit time). We indicate by $a^{\rm mp}, \theta^{\rm mp}, b^{\rm mp}$ the corresponding
quantities for a system with maximum packing.

Denote by $\mumin$ and $\mumax$ the vectors where all entries of $\mu$ are replaced by $\muminval
= \min_k \mu_k$ and $\mumaxval = \max_k \mu_k$, respectively, and let $r_\mu =
\mumaxval/\muminval$. Moreover, let us denote by $C_{m,n}$ the set of server configurations where
all layer-2 servers have been replaced by servers in layer~1, so that
\[
 C_{m,n} = \{ m' \in \Z_+^K: m'_k \ge m_k \ \forall k \ \ \text{and} \ \sum_k m'_k = \sum_k m_k + n  \}.
\]

\begin{theorem}
 \label{the:overallSteadyState}
 The stationary total number of customers in the system satisfies
 \begin{equation}
  \label{eq:overallQueueLength}
  |\bar X(m',0,\lambda,\mumax)| \lest |\bar X| \lest |\bar X^{\rm mp}(m,n,\lambda,\mumin)|
 \end{equation}
 for all $m' \in C_{m,n}$.
 Especially, the stationary mean number of customers is bounded by
 \begin{equation}
  \label{eq:overallWorkThroughput}
  \max_{m' \in C_{m,n}} a(m',0,\lambda,\mumax) \le \, a \, \le a^{\rm mp}(m,n,\lambda,\mumin),
 \end{equation}
 the mean throughput by
 \begin{equation}
  \label{eq:overallThroughput}
  \max_{m' \in C_{m,n}} r_\mu^{-1} \theta(m',0,\lambda,\mumax)
  \le \, \theta \,
  \le r_\mu \thetamp(m,n,\lambda,\mumin),
 \end{equation}
 and the overall blocking probability by
 \begin{multline}
  \label{eq:overallBlocking}
  1 - r_\mu(1-\bmp(m,n,\lambda,\mumin)) \\
  \le \, b \, \le \min_{m' \in C_{m,n}} \left( 1 - r_\mu^{-1} (1 - b(m',0,\lambda,\mumax)) \right).
 \end{multline}
\end{theorem}

\begin{remark}
 In the case where all service rates $\mu_k$ are equal, the bounds~\reff{eq:overallThroughput} and~\reff{eq:overallBlocking}
 can be written in a more natural form as
 \begin{alignat*}{2}
  \max_{m' \in C_{m,n}} \theta(m',0,\lambda,\mu)
  & \le \, \theta \,
  & & \le \thetamp(m,n,\lambda,\mu), \\
  \bmp(m,n,\lambda,\mu)
  &\le \, b \,
  & & \le \min_{m' \in C_{m,n}} b(m',0,\lambda,\mu).
 \end{alignat*}
\end{remark}

\begin{remark}
 The upper and lower bounds in~\reff{eq:overallQueueLength}, and hence the also the bounds in~\reff{eq:overallWorkThroughput}
 --~\reff{eq:overallBlocking}, are easy to compute numerically. The fast computation of the upper bound is explained
 in Remark~\ref{rem:maximumPacking}. To compute the
 lower bound, observe that $|\bar X(m',0,\lambda,\mumax)|$ has the same distribution as $\sum_k \erl(m'_k,\lambda_k/\mumaxval)$,
 where the terms in the sum are independent.
\end{remark}

\begin{proof}[Proof of Theorem~\ref{the:overallSteadyState}]
 Let $X$ be the number of customers in the original system, let $W$ be the
 number of customers in the system corresponding to the parameters $(m',0,\lambda,\mumax)$,
 and let $Y$ be the number of customers in the maximum packing system with parameters $(m,n,\lambda,\mumin)$.
 Assume that all processes are started at zero initial state. Then Theorem~\ref{the:packing} and
 Theorem~\ref{the:training} combined with Theorem~\ref{the:monotonemu} imply that
 \begin{equation}
  \label{eq:queuelengthTime}
  |W(t)| \lest |X(t)| \lest |Y(t)|
 \end{equation}
 for all $t$. Because all of the above processes are irreducible and positive recurrent, and because stochastic ordering is closed with respect to convergence in distribution~\cite{kamae1977}, taking $t\to \infty$ in~\reff{eq:queuelengthTime} shows the validity of~\reff{eq:overallQueueLength}. The bounds in~\reff{eq:overallWorkThroughput} follow by taking
 expectations, and the bounds in~\reff{eq:overallThroughput} from $\theta = \sum_k \mu_k a_k$. These bounds in
 turn imply~\reff{eq:overallBlocking}, because of the conservation law $(\sum_k \lambda_k)(1-b) = \theta$.
\end{proof}

\begin{figure}[h]
 \centering
 \psfrag{lambdanet}{\scriptsize $\lambdanet$}
 \psfrag{T}{\scriptsize $a$}
 \includegraphics[width=.49\textwidth]{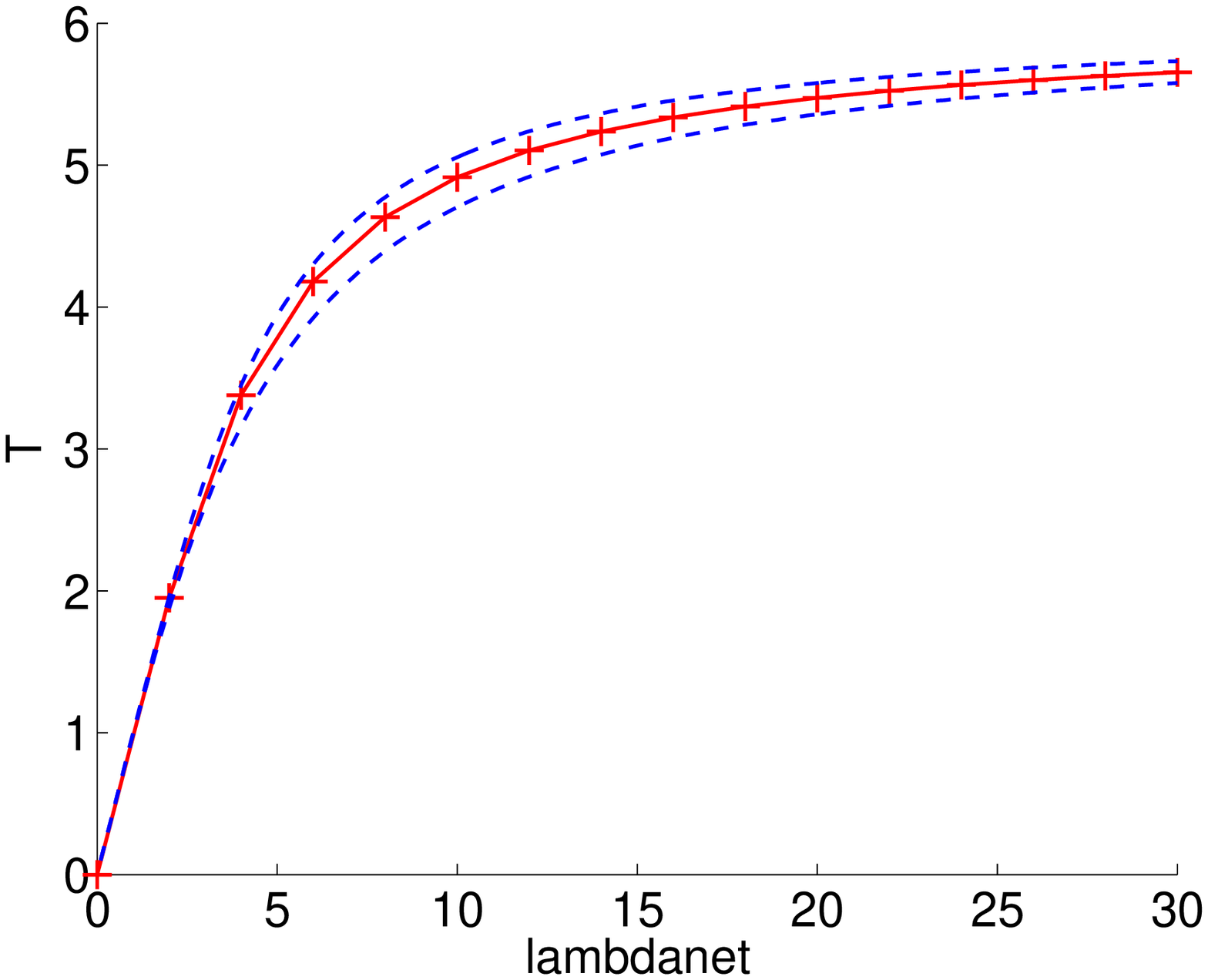}
 \includegraphics[width=.49\textwidth]{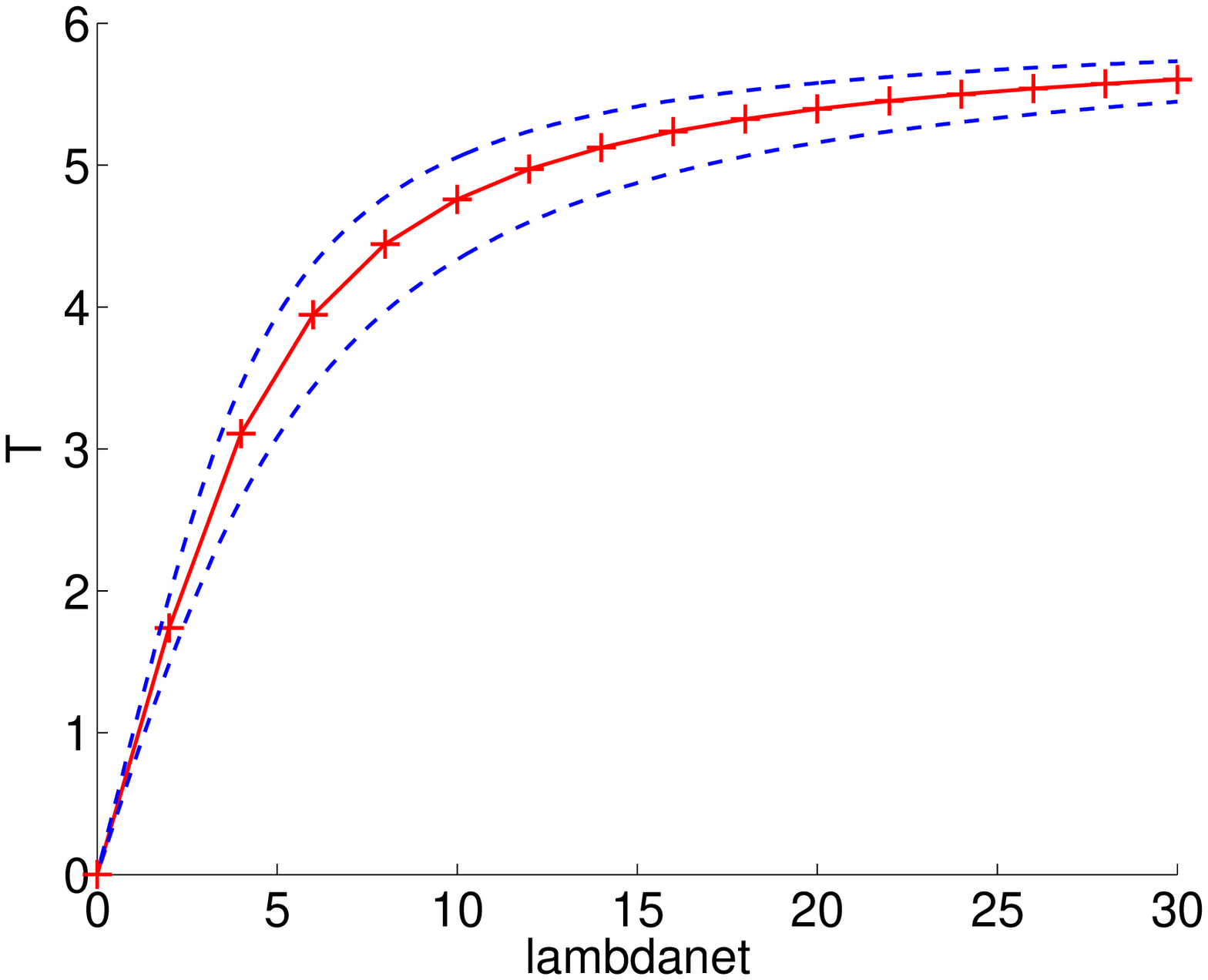}
 \caption{\label{fig:meanNumberOfCustomers} Stationary mean number of customers for $m=(2,2)$, $n=2$, and $\lambda = (\lambdanet/2,\lambdanet/2)$;
 where $\mu = (1,1)$ on the left, and $\mu = (1,1.3)$ on the right; plotted for the original system (solid line) and for the
 bounds in~\reff{eq:overallWorkThroughput} (dotted lines).}
\end{figure}
\begin{figure}[h]
 \centering
 \psfrag{ratioMu}{\scriptsize $\mu_2$}
 \psfrag{a}{\scriptsize $a$}
 \includegraphics[width=.49\textwidth]{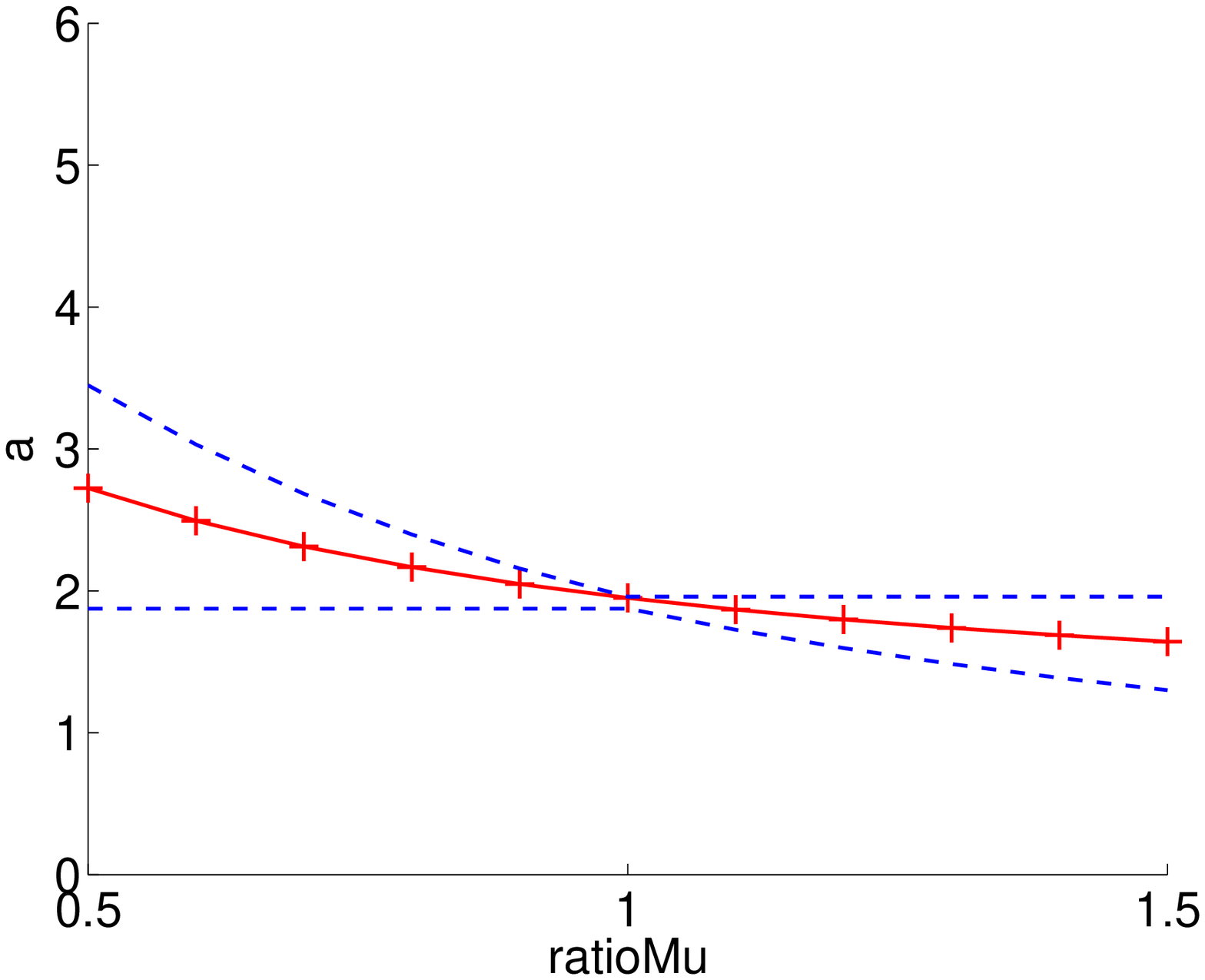}
 \includegraphics[width=.49\textwidth]{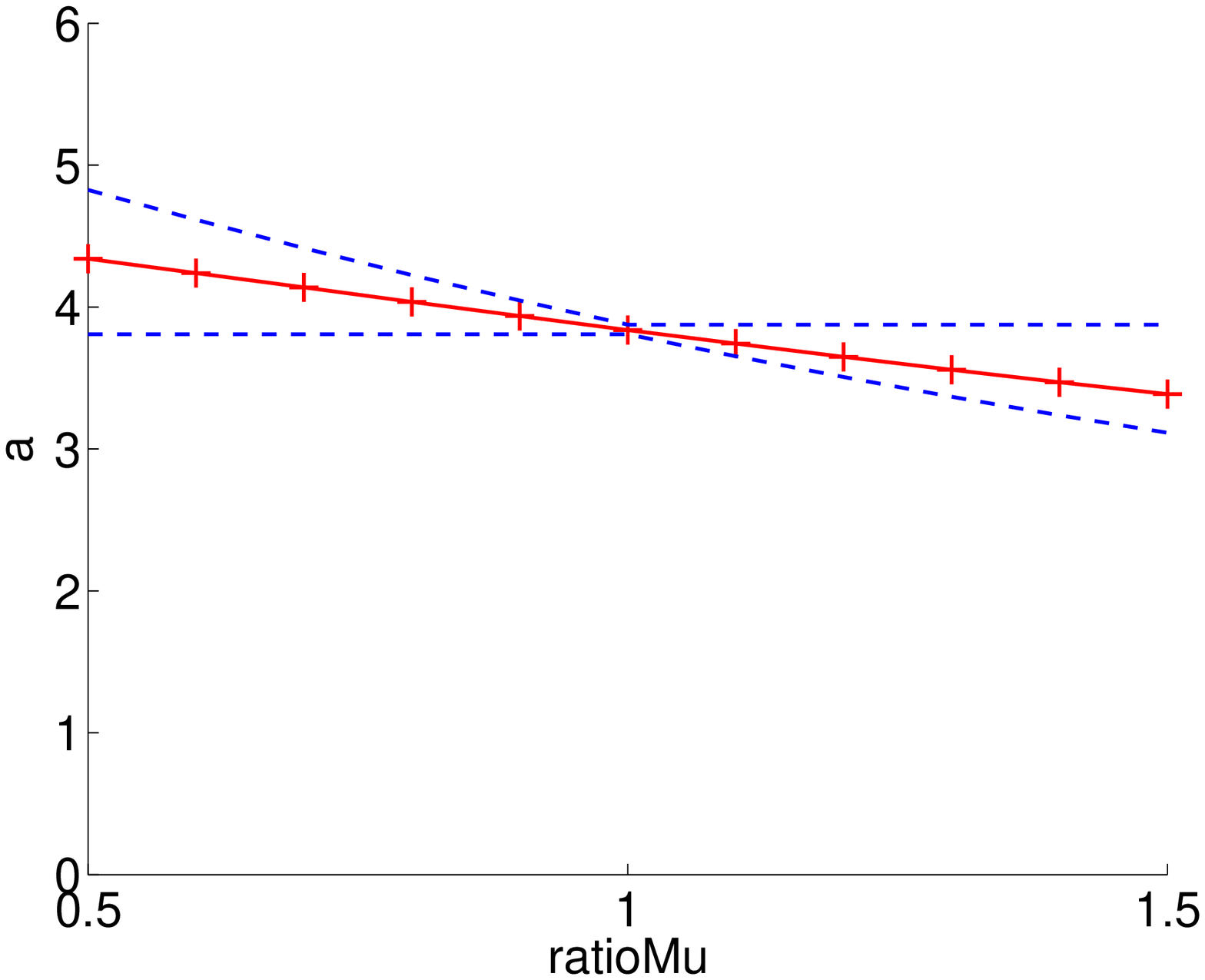}
 \caption{\label{fig:meanNumberOfCustomersVaryingMu} Stationary mean number of customers for $m=(2,2)$,
 $n=2$, $\mu_1 = 1$, and varying $\mu_2$, where $\lambda = (1,1)$ on the left and $\lambda = (1,5)$ on the right;
 plotted for the original system (solid line) and for the bounds in~\reff{eq:overallWorkThroughput} (dotted lines).}
\end{figure}
Figure~\ref{fig:meanNumberOfCustomers} illustrates the bounds~\reff{eq:overallWorkThroughput} of
the mean number of customers in a two-class system where the net arrival rate $\lambdanet$ is
varying. We see that the bounds are rather robust with respect to different values of the arrival
rates. Figure~\ref{fig:meanNumberOfCustomersVaryingMu} illustrates the same bounds for varying
$\mu_2$, showing that the accuracy of the bounds degrades rapidly as the difference of $\mu_2$ and
$\mu_1$ grows. This loss of accuracy is an inevitable consequence of replacing $\mu$ by $\mumin$
and $\mumax$ in~\reff{eq:overallWorkThroughput}. Intuitively one might think that the upper bounds
in~\reff{eq:overallQueueLength} and~\reff{eq:overallWorkThroughput} would hold without replacing
$\mu$ by $\mumin$. Example~\ref{exa:packing} shows that this is not true in general. However, the
right-hand side of~\reff{eq:overallWorkThroughput} with $\mu$ in place of $\mumin$, though not
generally an upper bound, appears to approximate well the original system for a wide range of
system parameters, even for the extreme choice of service rates of Example~\ref{exa:packing}. The
reason is that the actual repacking events in the maximum packing system occur relatively rarely
in moderately loaded systems; see Kelly~\cite{kelly1985} for an insightful discussion of this
phenomenon in the context of channel assignment in cellular radio networks.

\begin{example}
 \label{exa:packing}
 Consider a two-class loss network with server configuration $m=(1,0)$ and $n=2$. Assume $\lambda=(1,1)$ and
 $\mu = (\frac{1}{5},10)$, so that the service rates differ from each other by a factor of 50.
 Table~\ref{tab:packing} lists numerically calculated values of the stationary mean number of customers (per class and total)
 for the original loss network and the modification with maximum packing. The fact $a(m,n,\lambda,\mu) >
 \amp(m,n,\lambda,\mu)$ illustrates that for this special choice of parameters, maximum packing
 does \emph{not} increase the stationary mean number of customers in the system.
 \begin{table}[h]
  \centering
  \begin{tabular}[h]{l|r|r|r}
                                 & Class 1  & Class 2  & Total    \\ \hline
   $a(m,n,\lambda,\mu)$       & 2.325657 & 0.038612 & 2.364269 \\
   $\amp(m,n,\lambda,\mu)$    & 2.317818 & 0.046344 & 2.364162 \\
   $a(m,n,\lambda,\mumin)$    & 1.615744 & 0.997537 & 2.613281 \\
   $\amp(m,n,\lambda,\mumin)$ & 1.474617 & 1.172442 & 2.647059
  \end{tabular}
  \caption{\label{tab:packing} Mean number of customers in a loss network with and without maximum packing.}
 \end{table}
\end{example}

Example~\ref{exa:training} shows that replacing one layer-2 server by a layer-1 server may not decrease
the stationary mean number of customers, if not all service rates $\mu_k$ are equal. This shows that
it is necessary to replace $\mu$ by $\mumax$ in order to achieve a lower bound in~\reff{eq:overallQueueLength}.

\begin{example}
 \label{exa:training}
 Consider a two-class loss network with two different server configurations
 (i) $m=(0,0)$ and $n=3$, and (ii) $m'=(1,0)$, $n'=2$.
 Assume that $\lambda$ and $\mu$ are as in Example~\ref{exa:packing}.
 Numerically calculated values for the stationary mean number of customers (per class and total)
 given in Table~\ref{tab:training}.
 The fact $a(m,n,\lambda,\mu) < a(m',n',\lambda,\mu)$ illustrates that
 for this special choice of parameters, replacing one layer-2 server by a layer-1 server does \emph{not} decrease
 the stationary mean number of customers.
 \begin{table}[h]
  \centering
  \begin{tabular}[h]{l|r|r|r}
                                & Class 1  & Class 2  & Total    \\ \hline
   $a(m,n,\lambda,\mu)$      & 2.317808 & 0.046356 & 2.364164 \\
   $a(m',n',\lambda,\mu)$    & 2.325657 & 0.038612 & 2.364269 \\
   $a(m,n,\lambda,\mumax)$   & 0.099891 & 0.099891 & 0.199782 \\
   $a(m',n',\lambda,\mumax)$ & 0.099906 & 0.099453 & 0.199359
  \end{tabular}
  \caption{\label{tab:training} Mean number of customers in a loss network with two different server configurations.}
 \end{table}
\end{example}

\section{Conclusions}
\label{sec:conclusion}

Stochastic comparison techniques were developed for analyzing multiclass two-layer loss systems.
First, assuming all service rates to be equal, we proved that maximum packing stochastically
increases the total number of customers, and that moving a server from the second layer to the
first has the opposite effect. The monotonicity of the system with respect to service rates was
then used to extend the above conclusions to systems where the service rates may differ from each
other. As a consequence, computationally fast upper and lower bounds for the performance of the
system were derived.

The proofs of the main results (excluding Theorem~\ref{the:perClassLower}) were based on coupling
of continuous-time Markov processes, for which it was essential to assume that the service times
are exponentially distributed. On the other hand, the stationary distributions of the processes
acting as bounds in the main results, the maximum packing system and the Erlang loss system, are
known to be insensitive to the service time distribution~\cite{kelly1991}. This remarkable feature
calls for an extension of the comparison results to more general service time distributions. This
is an important open problem for which we believe that new probabilistic techniques are needed,
because a purely deterministic sample path approach was found unsuitable
(Examples~\ref{exa:noSamplePathRepacking} and~\ref{exa:noSamplePathPerClass}).

The accuracy of the bounds was numerically studied for systems with small number of servers. The
bounds for the per-class quantities appear not very accurate in general, though they may still be
useful in conservative dimensioning of system resources. The bounds for the aggregate system
quantities are much more accurate, especially when the mean service times across different
customer classes do not vary too much. For highly variable mean service times, the accuracy
degrades due to the need to modify the service time parameters in
Theorem~\ref{the:overallSteadyState}; however, if one uses the original $\mu$ in place of $\mumin$
in Theorem~\ref{the:overallSteadyState}, the maximum packing system appears to approximate well
the original system for a wide range of system parameters (though not anymore an upper bound in
general, see Example~\ref{exa:packing}). The accurate numerical evaluation of the system becomes
difficult when the number of servers is large, because of the rapid growth of the state
space~\cite{louth1994}. An interesting future problem is to asymptotically study the sharpness of
the bounds for large systems using scaling and renormalization techniques.

\section*{Acknowledgments}
We gratefully acknowledge helpful discussions with Sem Borst and sharp remarks by anonymous
referees. The main part of this research was carried out at Centrum voor Wiskunde en Informatica
and Eindhoven University of Technology. The research has been supported by the Dutch BSIK/BRICKS
PDC2.1 project, Helsingin Sanomat Foundation, and the Academy of Finland.

\bibliographystyle{apalike}
\bibliography{lslReferences}

\end{document}